\documentclass{amsart}
\usepackage{amsmath,amsfonts,amssymb,amsthm,amsxtra,latexsym,verbatim}

\newtheorem{theorem}{Theorem}[section]
\newtheorem{lemma}[theorem]{Lemma}
\newtheorem{crit}[theorem]{Criterion for Contribution}

\theoremstyle{definition}
\newtheorem{definition}[theorem]{Definition}
\newtheorem{example}[theorem]{Example}

\theoremstyle{remark}
\newtheorem{remark}[theorem]{Remark}

\numberwithin{equation}{section}

\newcommand{\A}{\mathbb{A}}
\renewcommand{\P}{\mathbb{P}}

\newcommand{\C}{\mathbb{C}}
\newcommand{\N}{\mathbb{N}}
\newcommand{\Q}{\mathbb{Q}}

\newcommand{\I}{\mathcal{I}}
\renewcommand{\O}{\mathcal{O}}

\newcommand{\J}{\mathfrak{J}}

\DeclareMathOperator{\ord}{ord}

\begin{document}

\title{Irrelevant Exceptional Divisors for Curves on a Smooth Surface}

\author{Karen E. Smith}
\address{Department of Mathematics, University of Michigan, Ann Arbor, Michigan 48109}
\email{kesmith@umich.edu}

\author{Howard M Thompson}
\address{Spring Arbor University, Spring Arbor, Michigan 49283}
\email{hthompso@arbor.edu}

\thanks{The authors were supported in part by NSF Grant DMS-0502170.}

\date{October 26, 2006 and, in revised form, July 26, 2007.}


\begin{abstract}
Given a singular curve on a smooth surface, we determine which exceptional divisors on the minimal resolution of that curve contribute toward its jumping numbers.
\end{abstract}

\maketitle

\section*{Introduction}

The \emph{jumping numbers} of a singular variety $Y$ embedded in a smooth complex variety $X$ form an interesting new invariant of the pair $(X, Y)$.  The jumping numbers are a sequence of positive rational numbers computed from--- and reflecting subtle information about--- an embedded resolution of singularities of the pair. For example, in the simplest case where $Y$ is a smooth hypersurface in $X$, the jumping numbers are simply the positive integers. But the sequence of jumping numbers becomes increasingly complicated as a resolution of singularities requires more blowings up or as functions vanishing on $Y$ vanish to higher orders along the resulting exceptional divisors. Jumping numbers, also known as jumping coefficients, were first explicitly defined in \cite{MR2068967} as those numbers $\lambda$ for which the multiplier ideal of the pair $(X, \lambda Y)$ makes a discrete "jump", though these natural invariants arose earlier in several contexts; see \cite{MR713242}, \cite{MR1056267}, and \cite{MR1198299}.

The jumping numbers of a pair $(X, Y)$ are determined by the exceptional divisors (or valuations) appearing in a log resolution.  On the other hand, some exceptional divisors never contribute to the jumping numbers. Since the jumping numbers do not depend on the choice of the log resolution, any divisor obtained by performing an "extraneous blowup" will be irrelevant from the point of view of jumping numbers. But examples show that even some "essential" exceptional divisors do not contribute to the sequence of jumping numbers. What is special about these divisors?

In this paper, we investigate this phenomenon, focusing on curves on a smooth surface, and establish precisely which exceptional divisors in a minimal log resolution of a singular curve on a smooth complex surface are "irrelevant" from the point of view of jumping numbers. Roughly stated, our main result is this: an exceptional divisor $E$ of a minimal embedded resolution of a curve $C$ contributes to the sequence of jumping numbers if and only if $E$ has non-trivial intersection with at least three of the (other) components of the full transform $C$. To state our result more precisely, we first review the definitions and make precise (see Definition~\ref{contribute}) the notion of an "exceptional divisor contributing a jumping number" to a pair. Our result can also be interpreted as saying that the divisors essential for computing jumping numbers are exactly those which appear on the log canonical model; see the last section of this paper. This interpretation is significant because it makes sense in any dimension (assuming the existence of log canonical models). In a forthcoming paper, we will return to this question in higher dimension using different techniques.

Our main theorem is very similar to results of of Favre and Jonsson~\cite{MR2138140}. Using completely different techniques, they give a valuative criterion for when a holomorphic function belongs to a multiplier ideal and observe that it is sufficient to check the valuations that correspond to Puiseux pairs. See Proposition~2.4, Lemma~2.11 and Fact~2 inside the proof of Theorem~6.1.

A variant of our question is asked by Koll\'ar in \cite[Section 10]{MR1492525}: which exceptional divisors contribute to the roots of the Bernstein-Sato polynomial of a complex polynomial?  Indeed, if $\lambda \leq 1$ is a jumping number of a polynomial $f$, then $-\lambda$ is a root of the Bernstein-Sato polynomial of $f$, as is shown in \cite{MR2068967}. Related questions have also been considered by Wim Veys, who considered which exceptional divisors contribute to Igusa (and motivic) Zeta functions.  We are grateful to Wim Veys for raising this issue to us in the context of jumping numbers.  We are also grateful to Jonathan Wahl, who brought to our attention the possible connection with log canonical models, and also to both Mark Spivakovsky and Bernard Teissier, with whom we consulted in understanding the connection with valuation theory described in Remark~\ref{val_rmk}. Also, we appreciated the remarks and support of Rob Lazarsfeld and J\'anos Koll\'ar as the project developed.

\section{Preliminaries.}
 
Let $Y$ be a proper closed subscheme of a smooth complex variety $X$, and let $\pi:\widetilde{X} \rightarrow X$ be a log resolution of the pair $(X, Y)$.  This means that $\widetilde X$ is smooth, that $\pi$ is a proper, birational morphism defining an isomorphism outside $Y$, and the ideal sheaf $\I_Y\O_{\widetilde X}$ is a locally principle ideal sheaf defining an effective divisor $G$ such that (together with the exceptional divisor of $\pi$) the reduced divisor of $G$ is in simple normal crossings.  Let $K_{\pi}$ or $K_{\widetilde{X}/X}$ denote the relative canonical divisor of the map $\pi$; this is the effective divisor on $\widetilde X$ defined by the Jacobian determinant of the morphism $\pi$.

\begin{definition}
Let $\lambda$ be any non-negative rational number.  The multiplier ideal of the pair $(X, Y)$ with coefficient $\lambda$ is defined as 
\[
\J(X, \lambda Y) = \pi_*\O_{\tilde X}(K_{\widetilde{X}/X} - \lfloor \lambda G \rfloor),
\]
where $\lfloor \lambda G \rfloor$ denotes the ``round-down'' of the $\Q$-divisor $\lambda G$, with $G$ defined as in
the preceding paragraph.
\end{definition}
\noindent
Note that, because both $G$ and $K_{\pi}$ are effective, 
\[
\J(X, \lambda Y) \subset \pi_*\O_{\tilde X}(K_{\widetilde{X}/X}) = \O_X,
\]
whence $\J(X, \lambda Y) $ really is an ideal sheaf on $X$. Furthermore, it is straightforward to check that it is independent of the choice of log resolution \cite[Theorem~9.2.18]{MR2095472}.
    
The following facts about multiplier ideals follow immediately from the definition:
\begin{itemize}
\item The multiplier ideals $\J(X, {\lambda}Y)$ are trivial for $\lambda$ sufficiently small;
\item If $\lambda_1 < \lambda_2$, then $\J(X, {\lambda_1}Y) \supseteq  \J(X, {\lambda_2}Y);$ 
\item For any $\epsilon > 0$ sufficiently small, $\J(X,{\lambda}Y) = \J(X, ({\lambda + \epsilon})Y )$ (because of the nature of rounding down);
\item For certain critical values of $\lambda$, we have a "jump" $\J(X, {\lambda}Y) \subsetneq \J(X, ({\lambda -\epsilon})Y)$.
\end{itemize}
This leads to the following definition.

\begin{definition}
A positive rational number $\lambda$ is a \emph{jumping number} of the pair $(X, Y)$ if $\J(X, {\lambda}Y) \subsetneq \J(X, ({\lambda - \epsilon})Y)$ for all $\epsilon > 0$.  The smallest jumping number is also called the \emph{log canonical threshold} of the pair $(X, Y)$.
\end{definition}

The jumping numbers of a pair $(X, Y)$ obviously form a discrete set of rational numbers. Indeed, if writing $G = \sum a_i E_i$ and $K_{\widetilde{X}/X} = \sum k_i E_i$, where the $E_i$ are effective prime divisors on $\tilde X$, then one easily verifies that the set of jumping numbers must be contained in the set of rational numbers of the form $\frac{k_i + n} {a_i}$, where $n$ is a positive integer.  Let us call the set of fractions of the form $\frac{k_i + n} {a_i}$ \emph{candidate jumping numbers because of} $E_i$. It is not true that every candidate jumping number is a jumping number, even for exceptional divisors in the minimal resolution of a plane curve singularity. The following example illustrates this.
   
\begin{example}
Let $C$ be the plane curve defined by the polynomial $x^4 - y^3 = 0$. The minimal log resolution $\pi: X \rightarrow \A^2$ has four exceptional divisors: $E_0$, obtained by blowing up the singular point\footnote{Here and throughout, we abuse notation by using the same notation to denote a divisor and its proper transform on any birational model.}; $E_1$, obtained by blowing up the intersection of $E_0$ with the curve $C$ (a point of second order tangency); $E_2$, obtained by blowing up the intersection of the three smooth curves $C$, $E_0$ and $E_1$ ($C$ is now a first order tangent to $E_0$ at this point); and finally $E_3$, obtained by blowing up the intersection of the three smooth curves $C$, $E_0$ and $E_2$. One computes
 \[
 \pi^*C = C + 3E_0 + 4E_1 + 8E_2 + 12E_3 \quad {\text{ and }}\quad  K_{\pi} = E_0 + 2E_1 + 4E_2 + 6E_3,
 \]
 so that the jumping numbers (say, those less than 1)  must all be of the form  
 \begin{itemize}
 \item  $\frac{1 + n}{3}$,  candidate jumping numbers from $E_0$;
 \item $\frac{2 + n}{4},$ candidate jumping numbers from $E_1$; 
 \item $\frac{4 + n}{8},$ candidate jumping numbers from $E_2$; or
 \item $ \frac{6 + n}{12}, $ candidate jumping numbers from $E_3$.
 \end{itemize}
On the other hand, direct computation shows that only the candidate jumping numbers from $E_3$ correspond to actual jumping numbers. In particular, $\frac{5}{8}$ and $\frac{7}{8}$ are not jumping numbers.
\end{example}
    
\begin{remark}
It is not hard to see that the jumping numbers of a pair $(X, Y)$ are periodic; this fact is essentially Skoda's theorem \cite[Theorem~9.3.24]{MR2095472}. In particular, when $Y$ is a divisor, the jumping numbers of $(X, Y)$ are completely determined by those jumping numbers less than one. Indeed, for any $\lambda >1 $, $\lambda$ is a jumping number if and only if $\lambda - 1$ is a jumping number whereas $1$ itself is always a jumping number. Therefore, in examining the jumping numbers of curves on surfaces, we may focus attention on jumping numbers strictly less than 1.
\end{remark}

\section{Jumping Numbers Contributed by a Divisor}
 
In this section, we make precise the idea that "a divisor $E_i$ contributes a jumping number $\lambda$" to the pair $(X, Y)$. The idea should be that we have a "jump"
\[
\J(X, {\lambda}Y) \subsetneq \J(X,  {(\lambda - \epsilon)}Y)
\]
occuring because the coefficient $k_i - \lfloor\lambda a_i \rfloor$ of $E_i$ has changed in passing from $\lambda$ to ($\lambda - \epsilon$)---and not because the corresponding coefficient of some other $E_j$ has changed.

Fix a variety $X$. A prime divisor $E$ lying on some variety $\widetilde X$ mapping properly and birationally to $X$ will be called a \emph{divisor centered on $X$.} We say that $E$ is an \emph{exceptional divisor centered on $X$} if its generic point contracts to a higher codimension point on $X$. Now, for a particular pair $(X, Y),$ we can ask whether a particular divisor $E$ centered on $X$ contributes to its jumping numbers. We propose the following definition.

\begin{definition} \label{contribute} 
Let $E$ be any prime divisor centered on $X$, and let $\lambda$ be a candidate jumping number for $E$. We say that $E$ {\emph contributes the jumping number $\lambda$ to the pair} $(X, Y)$ if
\[
\J(X, {\lambda}Y) \subsetneq \pi_*\O_{\tilde X}( K_{\pi} - \lfloor \lambda G\rfloor + E), 
\]
where $\pi: \widetilde{X} \rightarrow X$ is a log resolution of the pair $(X, Y)$ on which $E$ appears as a divisor, and $G$ is the divisor defined by $\I_Y \O_{\widetilde X} = \O_{\widetilde X}(-G) $.  A divisor $E$ is said to be \emph{relevant} for the pair (X, Y) if it contributes some jumping number and irrelevant otherwise.
\end{definition}

Definition~\ref{contribute} is independent of the choice of log resolution of $(X, Y)$.  Since the coefficient of $E$ in $\lambda G$ is an integer, rounding down has no effect and it depends only on the valuation determined by the divisor $E$, which of course is the same on any birational model.

If the coefficients of the divisors $E_i$ in the expression $G = \sum a_i E_i$ are all relatively prime, then for any particular $\lambda$, at most one of the coefficients in the expression $K_{\pi} - \lfloor \lambda G \rfloor$ can change as we pass from $\lambda$ to $\lambda - \epsilon$. In this case, the notion that ''the divisor $E_i$ contributes the jumping number $\lambda$" as in Definition~\ref{contribute} agrees with the intuition described immediately preceding it. On the other hand, some care is required when the $a_i$ are not relatively prime, as the next example shows.

\begin{example}\label{ex1}
Let $C$ be the plane curve defined by the polynomial $(x^3 - y^2)(x^2-y^3) = 0$, considered in a neighborhood of the origin.The minimal log resolution $\pi: X \rightarrow \A^2$ has five exceptional divisors: $E_0$, obtained by blowing up the singular point, $E_1$ and $E'_1$, obtained by blowing up the two intersections of $E_0$ with the curve $C$ (both points of tangency); and $E_2$ (respectively, $E'_2$), obtained by blowing up the intersection of the three smooth curves $C$, $E_0$ and $E_1$ (respectively, the three smooth curves $C$, $E'_0$ and $E'_1$). One computes
\begin{multline*}
\pi^*C = C + 4E_0 + 5(E_1+E'_1) + 10(E_2+E'_2) \\
\text{ and } K_{\pi} = E_0 + 2(E_1+E'_1) + 4(E_2+E'_2).
\end{multline*}
Here none of the exceptional divisors contributes the log canonical threshold, $\frac{1}{2}$ in the sense of Definition~\ref{contribute}.
\end{example}

In particular, rational numbers can be jumping numbers of a pair without being jumping numbers \emph{contributed by} any particular exceptional divisor. In Example~\ref{ex1}, the log canonical threshold is rather contributed by the reducible curve $E_0 + E_2 + E_2'$.  This issue is considered again in Section 4.

\section{Essential Exceptional Divisors for Plane curves}

First, a convention that will simplify notation. Let $C = \sum m_i C_i$ be a divisor on a smooth variety $X$, where each $C_i$ is a prime divisor.  Let $\pi: \widetilde{X} \rightarrow X$ be any proper birational morphism.  The proper transform of $C$ is the divisor $\sum m_i \widetilde C_i$ where $\widetilde C_i$ is the proper tranform on $X$ of the divisor $C_i$, that is, the closure of the image of the generic point of $C_i$ under the inverse birational map $\pi^{-1}$.  We will (abusively, but conveniently) use the notation $C$ to denote both the divisor on $X$ and its proper transform under any proper birational map.  

This section is devoted to the proof of the following theorem.
\begin{theorem} \label{main}
Let $C$ be a (possibly reducible and non-reduced) curve on a smooth complex surface $S$, and let $\pi: X \rightarrow S$ be the minimal embedded resolution of its singularities. Write
\[
\pi^* C =  C + \sum_i a_i E_i
\] 
where the $E_i$ are the exceptional divisors for $\pi$.  Then $E_j$ contributes jumping numbers to the pair $(S, C)$ if and only if
\[
(E_j \cdot E_j^{\circ}) \geq 3,
\]
where $E_j^{\circ}$ is the reduced divisor $ (\pi^*C)_{red} - E_j$.  In particular, if $(E_j\cdot E_j^{\circ}) \geq 3$, then $E_j$ contributes the jumping number $\lambda = 1 - 1/a_j$.
\end{theorem}

Fix any  exceptional divisor $E$. Consider the exact sequence
\[
0 \rightarrow \O_X(-E) \rightarrow \O_X \rightarrow \O_E \rightarrow 0. 
\]
Tensoring with the invertible sheaf $\O_X(K_{\pi} - \lfloor \lambda \pi^* C\rfloor + E)$ and taking advantage of the adjunction formula to compute $K_E$, we have an exact sequence
\[
0 \rightarrow \O_X(K_{\pi} - \lfloor \lambda \pi^* C\rfloor ) \rightarrow \O_X(K_{\pi} - \lfloor \lambda \pi^* C\rfloor + E) \rightarrow \O_E(K_{E} - \lfloor \lambda \pi^* C\rfloor \mid_{E}) \rightarrow 0.
\]
Pushing down to $S$, we arrive at
\[
0 \rightarrow \J(S, {\lambda}C) \hookrightarrow \pi_*\O_X (K_{\pi} - \lfloor \lambda \pi^*C \rfloor + E) \rightarrow \pi_*\O_E(K_E - \lfloor \lambda \pi^*C \rfloor\mid_{E}) \rightarrow 0,
\]
with the exactness on the right coming from the vanishing of \newline$R^1\pi_*\O_X(K_X - \lfloor \lambda\pi^*C \rfloor)$, essentially a consequence of Kawamata-Viehweg vanishing (see \cite[Theorem~9.4.1]{MR2095472}). Thus we see that $E$ contributes a jumping number if and only if
\[
\pi_*\O_E(K_E - \lfloor \lambda \pi^*C\rfloor \mid_{E}) = H^0(E, K_E - \lfloor \lambda \pi^*C) \rfloor \mid_{E}) 
\]
is non-zero, or equivalently, since $E = \P^1$, if and only if
\[
\deg( K_E - \lfloor \lambda \pi^*C\rfloor \mid_{E}) \geq 0.  
\]
Now, because  $K_{\P^1}$ has degree $-2$, we have proven the following
\begin{crit}\label{crit}
With notation as in Theorem~\ref{main}, the prime divisor $E$ contributes the jumping number $\lambda$ if and only if
\[
  -  \lfloor \lambda \pi^*C \rfloor \cdot E \geq 2.
\]
\end{crit}

To prove Theorem~\ref{main}, we also need the following Lemma.

\begin{lemma}\label{bound}
With notation as in Theorem~\ref{main}, fix $\lambda$ so that $\lambda a_j $ is an integer.  Then
\[
-  \lfloor \lambda \pi^* C \rfloor \cdot E_j   <  E_j^{\circ} \cdot E_j.
\]
\end{lemma}

\begin{proof} 
Let $C = \sum m_iC_i$.  We compute:
\begin{align*}
\lfloor\lambda\pi^*C\rfloor\cdot E_j 
&= \left(\sum\lfloor\lambda m_i \rfloor C_i + \sum\lfloor\lambda a_i \rfloor E_i \right)\cdot E_j \\
&= \left(\lambda a_j E_j + \sum\lfloor\lambda m_i \rfloor C_i + \sum_{i \neq j} \lfloor\lambda a_i \rfloor E_i \right)\cdot E_j\\
&> \left(\lambda a_j E_j + \sum (\lambda m_i - 1) C_i + \sum_{i \neq j} (\lambda a_i - 1) E_i \right)\cdot E_j\\
&= \lambda\left(\sum m_i C_i + \sum a_i E_i \right)\cdot E_j - \left(\sum C_i + \sum_{i \neq j} E_i \right)\cdot E_j \\
&= -\left(\sum C_i + \sum_{i \neq j} E_i \right)\cdot E_j.
\end{align*}

Here we have used that $\lambda a_j \in \N$ to get the second equality, that $\lfloor \gamma \rfloor > \gamma - 1$ for any rational number $\gamma$ to get the inequality, and that $\pi^*C \cdot E_j = 0$ to get the final equality.  The lemma is proved.
\end{proof}
   
We can now immediately deduce one implication of the theorem.  Suppose that $E_j^{\circ} \cdot E_j \leq 2$. From Lemma~\ref{bound}, we see that then $$- \lfloor \lambda \pi^* C \rfloor \cdot E_j < E_j^{\circ} \cdot E_j \leq 2,$$ for all candidate jumping numbers $\lambda$, so that $E_j$ does not contribute any jumping numbers by Criterion~\ref{crit}. This completes one direction of the argument.
 
For the other direction, assume that $E_j \cdot E_j^{\circ} \geq 3$. We will show then $E_j$ contributes the jumping number $\lambda = 1 - \frac{1} {a_j}$ to the pair $(S, C)$. We compute
\begin{align*}
\lfloor\lambda\pi^*C\rfloor\cdot E_j
&= \left(\sum\lfloor\lambda m_i \rfloor C_i + \sum\lfloor\lambda a_i \rfloor E_i \right)\cdot E_j \\
&= \left(\sum\left\lfloor m_i - \frac{m_i}{a_j} \right\rfloor C_i + \sum\left\lfloor a_i - \frac{a_i}{a_j}\right \rfloor E_i \right)\cdot E_j \\
&= \left(\sum m_i C_i + \sum a_i E_i \right)\cdot E_j + \left(\sum\left\lfloor -\frac{m_i}{a_j} \right\rfloor C_i + \sum\left\lfloor -\frac{a_i}{a_j} \right\rfloor E_i \right)\cdot E_j \\
&= -\left(\sum\left\lceil\frac{m_i}{a_j} \right\rceil C_i + \sum\left\lceil\frac{a_i}{a_j} \right\rceil E_i \right)\cdot E_j.
\end{align*}
Notice that whenever $C_i \cdot E_j \neq0$, $m_i \leq a_j$. So, for this choice of $\lambda$ we have that 
\[
  -    \lfloor\lambda\pi^*C\rfloor\cdot E_j =  \left(\sum C_i + \sum\left\lceil\frac{a_i}{a_j} \right\rceil E_i \right)\cdot E_j.
\]
Combining this computation with Criterion~\ref{crit}, we conclude that to show that $E_j$ contributes the jumping number $\lambda = 1 - \frac{1}{a_j}$ it is sufficient (and necessary) to show that
\[
\left(\sum C_i + \sum_i\left\lceil\frac{a_i}{a_j} \right\rceil E_i \right)\cdot E_j \geq 2.
\]
Under the hypothesis of Theorem~\ref{main} that $\left(\sum C_i + \sum_{i \neq j} E_i \right)\cdot E_j \geq 3$, this inequality evidently follows from
\[
-1 + \sum_{i \neq j}E_i \cdot E_j \leq \sum_{i}\left\lceil\frac{a_i}{a_j} \right\rceil E_i \cdot E_j,
\]
which in turn is equivalent to 
\begin{equation}\label{e2}
  -1 - E_j^2 \leq \sum_{i\neq j}\left(\left\lceil\frac{a_i}{a_j}  \right\rceil - 1\right)  E_i \cdot E_j.
\end{equation}
We will be able to verify inequality~\eqref{e2} after we have computed $E_j^2$.

We now compute $E_j^2$.  Imagine the minimal log resolution $\pi: X \rightarrow S$ as a sequence of point blowups.  We say that $E_i$ is proximate to $E_j$ or $E_i$ \emph{is created from} $E_j$ if, in the process of carrying out these point blowups to get $\pi$, the divisor $E_i$ is created by blowing up a point on $E_j$. The notion of  proximity was first introduced by Enriques and Chisini (see for example \cite{MR1266187}). Lemma~\ref{lem} and Lemma~\ref{lem11} are well known in the context of proximity. Their proofs are included to make the proof of Theorem~\ref{main} self-contained. Note that this relationship does not imply that $E_i$ and $E_j$ intersect on the log resolution $X$, because we may also have to blow up their intersection point.
  
\begin{lemma}\label{lem}
Let $E$ be any prime exceptional divisor for a proper birational map of smooth surfaces $\pi: X \rightarrow S$, and let $n$ be the number of exceptional divisors created from $E$. Then
\[
E^2  = -n -1. 
\]
\end{lemma}

\begin{proof}
The proof is an easy exercise using induction on $n$. The inductive step is accomplished using the following observation: if $f: X_2 \rightarrow X_1$ is the blowup of a point on smooth divisor $E$ on $X_1$, and $F$ is the corresponding exceptional divisor, then
\[
E^2 \,\,{\text{ on }}\,\,X_1 = (f^* E)^2 \,\,{\text{ on }}\,\, X_2 = (F + E)^2 = F^2 + 2F\cdot E + E^2 = E^2 + 1 \,\,\, {\text{ on }}\,\, X_2.
\]
\end{proof}

Using Lemma~\ref{lem},  we see that inequality~\eqref{e2} is equivalent to 
\begin{equation}\label{e3}
n \leq \sum_{i\neq j}\left(\left\lceil\frac{a_i}{a_j} \right\rceil - 1\right)  E_i \cdot E_j, 
\end{equation}
where $n$ denotes the number of exceptional divisors created created from $E_j$ in the minimal log resolution $\pi$. Since
\[
\left\lceil\frac{a_i}{a_j}  \right\rceil - 1 \geq 0
\]
in any case, it suffices to show that for each of the $n$ divisors $E_i$ created from $E_j$, we get some \emph{positive} contribution to the summation in \eqref{e3} above.  If all such $E_i$ would have non-zero intersection with $E_j$, it would suffice to show that $a_i > a_j$ whenever $E_i$ is created from $E_j$, for then each of the $n$ coefficients $\lceil\frac{a_i}{a_j} \rceil - 1$ would be strictly positive. Unfortunately, however, this is not the case: it can happen that $E_i \cdot E_j = 0$ even when $E_i$ is created by blowing up a point on $E_j$.  Indeed, this happens precisely when a further exceptional divisor is created by blowing up the intersection point of $E_i$ and $E_j$. In general, the set of exceptional divisors of $\pi$ created from $E_j$ can be partitioned into chains of exceptional divisors (each created by blowing up a point on its predecessor) in which only the last divisor in the chain intersects $E_j$.  It therefore suffices to prove the following Lemma.
      
\begin{lemma}\label{lem11}
Let $E_1, \dots, E_s$ be a connected chain of exceptional divisors on the minimal log resolution $\pi$ of $(S, C)$, all obtained by blowing up points on $E_j$. Then, for all $i$, $1 \leq i \leq s$,
\[
\ord_{E_i}(\pi^* C) > i  \,  \ord_{E_j}(\pi^* C),
\]
where we have ordered the divisors so that $ E_{i+1}$ is created by blowing up the intersection point $E_i \cap E_j$.
\end{lemma}
      
\begin{proof}
The result easily follows by induction after one observes the following basic facts.  First, once a divisor is created, its multiplicity in the full transform of $C$ is unchanged by further blowing up, so to understand $\ord_{E_i}(\pi^* C)$, we can focus attention on the blowup $f: X_2 \rightarrow X_1$ that creates $E_i$. Second, when we create $E_i$ by blowing up the intersection point $p$ of $E_{j}$ and some other prime divisor $D$ in the full transform of $C$ on $X_1$, we always have that
\[
\ord_{E_i}(f^*g^* C) \geq \ord_{E_j}(g^*C) + \ord_D(g^*C),
\]
where $g$ denotes the map $X_1 \rightarrow X$.  (Of course, the multiplicity of $E_i$ in the full transform $f^*g^*C$ of $C$ on $X_2$ will be larger than this sum if the point $p$ also lies on some other divisor in $g^*C$ or if $p$ is a singular point on $D$.)  Finally, because $\pi$ is a \emph{minimal} log resolution of $(S, C)$, we only blow up points $p$ that lie on $C$.

Putting these observations together, we see that first
\[
\ord_{E_1 }(\pi^*C) \geq  \ord_{E_j}(\pi^*C) + \ord_{C_k}(\pi^*C) > a_j,
\]
where $E_1$ is obtained by blowing up $p$ which lies on $E_j$ and some component $C_k$ of $C$. Assume inductively that $ \ord_{E_{i-1} }(\pi^*C) > (i-1) \ord_{E_{j} }(\pi^*C)$.  Then blowing up the intersection point $E_{i-1}  \cap E_j$ to create $E_i$, we have
\begin{align*} 
\ord_{E_i}(\pi^*C) &\geq \ord_{E_{i-1} }(\pi^*C) + \ord_{E_{j} }(\pi^*C) \\
&> (i-1) \ord_{E_{j} }(\pi^*C) + \ord_{E_{j} }(\pi^*C) \\
&= (i) \ord_{E_{j} }(\pi^*C).
\end{align*}
The proof of the lemma, and hence Theorem~\ref{main},  is complete.
\end{proof}

\begin{remark}\label{val_rmk}
It is also possible to describe the relevant divisors for computing jumping numbers in terms of the Puiseux data of a curve. For simplicity, we assume that $(S, C)$ is a plane branch, which is to say, the germ of an analytically irreducible curve in the plane. Let $\O_C = \C[[x, y]]/(f)$ be the complete local ring of $C$ at the singular point, and consider the normalization map $\O_C \hookrightarrow \C[[t]]$, where $t$ is some uniformization parameter.  This map determines an obvious valuation $\nu_g$ on $\C((x, y))$, whose value on a power series $h$ is simply the order of $t$ in the image of $h$ under the composition $\C[[x, y]] \rightarrow \O_C \hookrightarrow \C[[t]]$.
  
Let $\Gamma = \langle \beta_0, \beta_1, \dots, \beta_g \rangle$ be minimal generators of the sub-semi-group of $\N$ generated by the values under $\nu_g$ of the elements in $\C[[x, y]]$.  Choosing coordinates appropriately, we may assume $\nu_g(x) = \beta_0,\,\, \nu_g(y) = \beta_1 > \beta_0$, and we may fix $f_i \in \C[[x, y]]$ such that $\nu_g(f_i) = \beta_i$ for $0 \leq i \leq g$ (by convention $f_1 = y$ and $f_0 = x$).  This data determines $g$ valuations on $\C((x, y))$, successively dominating each other: $\nu_1$ is the monomial valuation taking the values $\beta_0$ and $\beta_1$ on $x$ and $y$ respectively, while $\nu_k$  is the valuation uniquely determined by the data that $\nu_k(f_i) = \beta_i$ for $0 \leq i \leq k$. In particular, $v_g$ agrees with the original valuation we constructed by normalization. See section eight of \cite{MR1037606}, especially Remark~8.14 on page 148.
 
The $g$-valuations we constructed this way correspond exactly to the valuations along the certain exceptional divisors $E_i$ appearing in a minimal log resolution of the singularities of $C$---- indeed, these are precisely the exceptional divisors $E$ having the property that $E \cdot E^{\circ} \leq 3$ \cite{MR1037606}.  Thus the divisors relevant from the Puiseux series point of view are precisely the divisors relevant for computing jumping numbers.

It is possible to give a precise formula for the jumping numbers of an analytically irreducible curve $(S, C)$ in terms of the data $\beta_0, \dots, \beta_g$. This is done in the forthcoming PhD thesis of Tarmo J\"arvilehto~\cite{TJ07}.
\end{remark}

\section{A speculation about the higher dimensional case.}

The divisors relevant for computing jumping numbers of curves on surfaces (Theorem~\ref{main}) have an alternative interpretation in the language of the minimal model program, which suggests an interpretation for relevant divisors in higher dimension as well.

Let $(S,C)$ be a pair consisting of a divisor on a smooth surface, and let $\pi: X \rightarrow S$ be a log resolution, with $\Delta$ denoting the reduced divisor of $\pi^*C$.  The pair $(X, \Delta)$ admits a unique $S$-relative log canonical model $(\tilde X, \tilde \Delta)$, which is independent of the original choice of log resolution of $(S, C)$; see \cite[Section 3.8]{MR1658959}.  The divisors relevant for jumping numbers of $(S, C)$ are \emph{precisely} those appearing on this $S$-relative log canonical model.

There are several ways to see this.  Assuming that $X$ is a \emph{minimal} log resolution of $(S, C)$, the morphism $\nu:X\to \tilde X$ collapses precisely those exceptional divisors $E$ on $X$ for which $(K_{X/S} + \Delta)\cdot E \leq 0$, or equivalently, precisely those divisors $E$ such that $E \cdot E^{\circ} \leq 2$ in the notation of Theorem~\ref{main} (See \cite[Prop~2.5]{MR1425327}).  Using this, one can argue that in the surface case, the divisors relevant to computing jumping numbers can be described as those divisors surviving on the log canonical model of $(X, (\pi^* C)_{red})$, where $\pi$ is any log resolution of $(S, C)$.
  
In a forthcoming paper, we will give a different argument valid in any dimension to prove that the multiplier ideal of a pair can be computed directly from any log minimal model of $(X, \Delta)$. This shows in particular that any divisors that are collapsed to higher codimension on some log minimal model can not contribute to the jumping numbers. We hope to show eventually that in fact the multiplier ideal can be computed from the log canonical model, and that indeed, every divisor appearing on the log canonical model contributes some jumping number. For now we can only conjecture that the divisors relevant for jumping numbers in general are precisely the divisors surviving on the log canonical model.

\bibliographystyle{amsalpha}


\begin{thebibliography}{ELSV04}

\bibitem[ELSV04]{MR2068967}
Lawrence Ein, Robert Lazarsfeld, Karen~E. Smith, and Dror Varolin,
  \emph{Jumping coefficients of multiplier ideals}, Duke Math. J. \textbf{123}
  (2004), no.~3, 469--506. \MR{MR2068967 (2005k:14004)}

\bibitem[FJ05]{MR2138140}
Charles Favre and Mattias Jonsson, \emph{Valuations and multiplier ideals}, J.
  Amer. Math. Soc. \textbf{18} (2005), no.~3, 655--684 (electronic).
  \MR{MR2138140 (2007b:14004)}

\bibitem[J\"07]{TJ07}
Tarmo J\"{a}rvilehto, \emph{Jumping numbers of a simple complete ideal in a
  two-dimensional regular local ring}, Ph.D. thesis, University of Helsinky,
  2007.

\bibitem[KM98]{MR1658959}
J{\'a}nos Koll{\'a}r and Shigefumi Mori, \emph{Birational geometry of algebraic
  varieties}, Cambridge Tracts in Mathematics, vol. 134, Cambridge University
  Press, Cambridge, 1998, With the collaboration of C. H. Clemens and A. Corti,
  Translated from the 1998 Japanese original. \MR{MR1658959 (2000b:14018)}

\bibitem[Kol97]{MR1492525}
J{\'a}nos Koll{\'a}r, \emph{Singularities of pairs}, Algebraic geometry---Santa
  Cruz 1995, Proc. Sympos. Pure Math., vol.~62, Amer. Math. Soc., Providence,
  RI, 1997, pp.~221--287. \MR{MR1492525 (99m:14033)}

\bibitem[Laz04]{MR2095472}
Robert Lazarsfeld, \emph{Positivity in algebraic geometry. {II}}, Ergebnisse
  der Mathematik und ihrer Grenzgebiete. 3. Folge. A Series of Modern Surveys
  in Mathematics [Results in Mathematics and Related Areas. 3rd Series. A
  Series of Modern Surveys in Mathematics], vol.~49, Springer-Verlag, Berlin,
  2004, Positivity for vector bundles, and multiplier ideals. \MR{MR2095472
  (2005k:14001b)}

\bibitem[Lib83]{MR713242}
A.~Libgober, \emph{Alexander invariants of plane algebraic curves},
  Singularities, Part 2 (Arcata, Calif., 1981), Proc. Sympos. Pure Math.,
  vol.~40, Amer. Math. Soc., Providence, RI, 1983, pp.~135--143. \MR{MR713242
  (85h:14017)}

\bibitem[Lip94]{MR1266187}
Joseph Lipman, \emph{Proximity inequalities for complete ideals in
  two-dimensional regular local rings}, Commutative algebra: syzygies,
  multiplicities, and birational algebra (South Hadley, MA, 1992), Contemp.
  Math., vol. 159, Amer. Math. Soc., Providence, RI, 1994, pp.~293--306.
  \MR{MR1266187 (95j:13018)}

\bibitem[LV90]{MR1056267}
F.~Loeser and M.~Vaqui{\'e}, \emph{Le polyn\^ome d'{A}lexander d'une courbe
  plane projective}, Topology \textbf{29} (1990), no.~2, 163--173.
  \MR{MR1056267 (91d:32053)}

\bibitem[Spi90]{MR1037606}
Mark Spivakovsky, \emph{Valuations in function fields of surfaces}, Amer. J.
  Math. \textbf{112} (1990), no.~1, 107--156. \MR{MR1037606 (91c:14037)}

\bibitem[Vaq92]{MR1198299}
Michel Vaqui{\'e}, \emph{Irr\'egularit\'e des rev\^etements cycliques des
  surfaces projectives non singuli\`eres}, Amer. J. Math. \textbf{114} (1992),
  no.~6, 1187--1199. \MR{MR1198299 (94d:14015)}

\bibitem[Vey97]{MR1425327}
Willem Veys, \emph{Zeta functions for curves and log canonical models}, Proc.
  London Math. Soc. (3) \textbf{74} (1997), no.~2, 360--378. \MR{MR1425327
  (98k:11176)}

\end{thebibliography}

\providecommand{\bysame}{\leavevmode\hbox to3em{\hrulefill}\thinspace}
\providecommand{\MR}{\relax\ifhmode\unskip\space\fi MR }
\providecommand{\MRhref}[2]{%
  \href{http://www.ams.org/mathscinet-getitem?mr=#1}{#2}
}
\providecommand{\href}[2]{#2}

\end{document}